\magnification=1200

\tolerance=500

\font \twelvebf=cmbx12

\hfuzz=3pt

\def \a {\alpha}

\def \avd {anneau de valuation discr\`ete}

\def \Coker {{\rm Coker\,}}

\def\dem{\noindent {\bf D\'emonstration.}\enspace \nobreak }

\def \expl#1 {\medbreak \noindent {\bf Exemple #1.}\enspace }
\def \expls#1 {\medbreak \noindent {\bf Exemples #1.}\enspace }

\def \Ext {{\rm Ext\,}}

\def \fl {\rightarrow }
\def \Fl#1{\smash{\mathop{\longrightarrow}\limits^{#1}}}

\def \Hom { {\rm Hom\,}}

\def \id #1 {<\!\! #1 \!\! >}
\def \Im {{ \rm Im\,}}

\def \impl {\Rightarrow}

\def \la {\lambda}

\def \lign {\hfil \break }

\def \Vf#1{\Big\downarrow \rlap{$\vcenter{\hbox{$\scriptstyle#1$}}$}}

\def \mud {{\,\mid\,}}

\def \ov {\overline}

\def \rema#1 {\medbreak \noindent {\bf Remarque #1.}\enspace }
\def \remas#1 {\medbreak \noindent {\bf Remarques #1.}\enspace }

\def \rg {{ \rm rang \,}}

\def \s {\sigma}

\def \sh {{^\sharp}}
\def \Spec {{\rm Spec\,}}

\def \T {\otimes }
\def \tarte#1 {\medbreak \noindent {\it #1.}\medbreak}

\def \titre#1{\medbreak \noindent {\bf #1.}\medbreak}

\outer \def \th #1. #2\par{ \medbreak 
\noindent {\bf#1. \enspace} {\sl#2 }\par
\ifdim \lastskip< \medskipamount \removelastskip \penalty55 \medskip \fi}

\def  \Vf#1{\Big\downarrow \rlap{$\vcenter{\hbox{$\scriptstyle#1$}}$}}
\def \vf {\downarrow}

\def \wi {\widetilde}

\def \cC {{\cal C}}
\def \cE {{\cal E}}
\def \cL {{\cal L}}
\def \cP {{\cal P}}
\def \cF {{\cal F}}
\def \cG {{\cal G}}
\def \cH {{\cal H}}
\def \cJ {{\cal J}}
\def \cO {{\cal O}}
\def \cK {{\cal K}}

\def \cN {{\cal N}}

\def \bP {{\bf P}}

\def \bZ {{\bf Z}}
\def \bN {{\bf N}}

\def \sN{{\sl N}}

\def \Om {{\Omega}}

\def \sN {{\sl N}}

\centerline {\twelvebf Construction de familles minimales de courbes gauches}

\vskip 2 cm

\titre { Introduction}

Soit $A$ un anneau local noeth\'erien, $t$ le point ferm\'e de $T = \Spec A$ et soit 
$\bP^3_A$ l'espace projectif associ\'e. Si $\cC$ est une famille plate de courbes de
$\bP^3_A$, il existe une r\'esolution  de $\cC$, dite
de type \sN, c'est-\`a-dire une suite exacte
$0 \fl \cP \fl \cN \fl \cJ_\cC \fl 0$ avec $\cN$ localement libre et $\cP$
dissoci\'e  (somme directe de faisceaux inversibles) et le faisceau $\cN$ est bien
d\'etermin\'e par $\cC$ \`a pseudo-isomorphisme pr\`es (cf. [HMDP1]), de
sorte qu'on a une application $\Phi$ qui \`a une courbe $\cC$ associe la classe de
pseudo-isomorphisme du faisceau $\cN$. Dans [HMDP1] nous avons d\'etermin\'e les
fibres de $\Phi$ \`a d\'ecalage pr\`es  en montrant  que les faisceaux
$\cN=
\Phi (\cC)$ et
$\cN'= \Phi (\cC')$  sont pseudo-isomorphes \`a d\'ecalage pr\`es si et seulement si  les
courbes
$\cC$ et
$\cC'$ sont dans la m\^eme classe de biliaison  (c'est une g\'en\'eralisation du
th\'eor\`eme de Rao, cf. [R]).

\vskip 0.3 cm

Par ailleurs, il est facile de voir que $\Phi$ est surjective {\it \`a d\'ecalage pr\`es}
(cf.  2.7 ci-dessous ou
[S1]), c'est-\`a-dire que  si $\cN$  est un faisceau localement libre sur $\bP^3_A$ il
existe une famille de courbes  (lisses et connexes)
$\cC$  de
$\bP^3_A$, plate  sur $A$, et un entier $h$ tels que $\cC$ admette une r\'esolution de
type {\sl N} de la forme :
$ 0 \fl \cP \fl \cN \fl \cJ_\cC (h) \fl 0.$

\vskip 0.3 cm

L'int\'er\^et de ce type de r\'esultat d'existence provient de l'\'etude  du sch\'ema de
Hilbert 
$H_{d,g}$ des courbes de degr\'e $d$ et genre $g$ donn\'es et notamment
 des sp\'ecialisations qu'on peut obtenir
 lorsqu'on
fait varier la cohomologie. Une question essentielle est de montrer l'existence d'une
famille de courbes, param\'etr\'ee par un anneau de valuation discr\`ete
$A$, dont le point sp\'ecial 
$C_0$ (resp. le point g\'en\'erique $C$) est une courbe \`a  cohomologie 
pres\-crite. La m\'ethode que nous proposons pour cela s'inspire de la construction des
courbes
\`a partir des modules de Rao inaugur\'ee dans [MDP1] IV. Elle consiste,  \`a partir d'un
objet (appel\'e triade) qui g\'en\'eralise le module de Rao, cf. [HMDP3],
\`a d\'eterminer un faisceau
$\cN$ convenable puis
\`a construire la famille de courbes avec une r\'esolution de type \sN\ comme
ci-dessus.

Pour ce type de probl\`emes,  le  r\'esultat de 2.7 est insuffisant, puisque, s'il assure
l'existence  d'une courbe de la classe de biliaison associ\'ee \`a  $\cN$, il ne permet
pas de contr\^oler le d\'ecalage
$h$ et donc d'affirmer que la famille construite a les  degr\'e et genre voulus.
Dans cet article nous d\'eterminons toutes les familles de courbes de la classe de
biliaison associ\'ee \`a un faisceau localement libre $\cN$ donn\'e, \`a commencer par
les familles minimales, i.e. celles pour lesquelles l'entier $h$ est minimal.

 Pour cela nous introduisons ci-dessous, cf. 2.4, une fonction
$q: \bZ \fl \bN$ et un entier $b_0$ associ\'es \`a ce faisceau. Cette fonction est
d\'efinie  \`a
la mani\`ere de [MDP1,3,4] comme une fonction $q$ du faisceau $\cN_t = \cN \T_A k(t)$, 
associ\'ee non pas au module des sections globales $H^0_* \cN_t$ tout entier, mais  au
sous-module des sections  qui pro\-vien\-nent de
$H^0_*\cN$, cf. \S1. Les familles de courbes  dont la r\'esolution de type {\sl N} 
correspond
\`a
$\cN$ sont alors donn\'ees par le th\'eor\`eme suivant : 

\th {Th\'eor\`eme 2.6}.   Soit $\cN$   un faisceau sur $\bP^3_A$
localement libre de rang $r$ et non dissoci\'e.  Soit
$\cP =
\bigoplus_{n
\in \bZ} \cO_{\bP_A} (-n)^{p(n)}$ un faisceau dissoci\'e  de rang $r-1$.  Alors,
les conditions suivantes sont \'equivalentes : \lign i) Pour $u$ g\'en\'eral
dans 
$ H=\Hom_{\cO_{\bP}} (\cP,
\cN)$, on a une suite exacte :
$$0 \fl \cP \Fl {u} \cN \fl \cJ_{\cC}(h) \fl 0 \leqno {(*)}$$
 o\`u $\cC$
est une famille de courbes de $\bP^3_A$, plate  sur $A$. \lign
ii) 
La fonction $p$ v\'erifie les conditions suivantes : \lign
\indent 1) on a $p\sh (n) \leq q\sh (n)$ pour tout $n \in \bZ$, \lign
\indent 2) s'il existe $n \leq b_0$ tel que l'on ait $p\sh (n)=q\sh (n)$  on a un
isomorphisme $(\cP_{\leq n})_t \simeq \cN_{t,\leq n}$ (cf. 2.4). \lign
La valeur minimale du d\'ecalage $h$ est \'egale
\`a
$\sum_{n\in
\bZ} nq(n) +
\deg
\cN$. La famille de courbes correspondant \`a ce d\'ecalage minimum est donn\'ee par une
r\'esolution comme ci-dessus, avec $p=q$.

\vskip 0.3 cm

Les deux  corollaires suivants d\'ecrivent les courbes  d'une classe de biliaison
donn\'ee :

\th {Corollaire 2.9}.  Soit $\cN_0$ un faisceau localement libre
extraverti (i.e. tel que $H^1_* \cN_0^\vee$  soit nul) minimal (i.e.  sans facteur direct
dissoci\'e). \lign
Posons   $q=
q_{\cN_0}$ et $h_0= \sum_{n \in \bZ} nq(n) + \deg \cN_0$. \lign
1)  Il existe une
famille de courbes  $\cC_0$ et une r\'esolution 
$ 0 \fl \cP_0 \Fl{v} \cN_0 \fl \cJ_{\cC_0} (h_0) \fl 0$ avec $\cP_0$ dissoci\'e. \lign
2) Si $\cC_1$ est une famille
de courbes  de la classe de biliaison de $\cC_0$, elle admet une r\'esolution 
$ 0 \fl \cP \fl \cN_0 \oplus \cL \fl \cJ_{\cC_1} (h) \fl 0$ avec $\cL$ dissoci\'e et on a
$h
\geq h_0$. Si $d$ et $g$ (resp. $d_0$ et $g_0$) sont respectivement le degr\'e et le
genre de $\cC_1$ (resp. $\cC_0$) on a $d \geq d_0$ et $g \geq g_0$.\lign
3) R\'eciproquement, pour tout $h \geq h_0$, il existe une famille de courbes $\cC_1$
avec une r\'esolution comme ci-dessus.
\lign
4) Si de plus on a
$h=h_0$, $\cC_0$ et
$\cC_1$ sont jointes par une d\'eformation  \`a cohomologie uniforme (cf.  2.8
ci-dessous). On dit que
$\cC_0$ est une {\bf famille minimale} de courbes.

\th {Corollaire 2.10}. (Propri\'et\'e de Lazarsfeld-Rao)  Soit
$\cC_0$ une famille minimale de courbes et $\cC$  une famille de courbes de la classe de
biliaison de
$\cC_0$. Alors, il existe un entier $m \geq 0$ et une suite de courbes $\cC_0, \cC_1,
\cdots, \cC_m$ telle que $\cC_{i+1}$ s'obtienne \`a partir de $\cC_i$ par une
biliaison \'el\'ementaire  (cf. [HMDP1] ) et $\cC$ \`a partir de $\cC_m$ par
une d\'eformation  \`a cohomologie uniforme.

\vskip 0.3 cm

 Lorsque le faisceau
$\cN$ est \'ecrit comme image d'un morphisme de faisceaux dissoci\'es $ s: \cL_2 \fl
\cL_1$ tel que la fl\`eche $L_2 \fl H^0_* \cN$ induite par $s$ soit surjective et que le
conoyau de $s$ soit localement libre (ce qui est   possible pour tout faisceau
localement libre sur $\bP^3_A$), nous donnons  un algorithme  de calcul de la fonction
$q$, enti\`erement effectif en termes  de certains mineurs de $s_t$, (cf. 3.1).
Nous donnons  quelques exemples  simples d'applications de cet algorithme.
Pour d'autres exemples et des applications de ces techniques dans le cadre des triades,
le lecteur se reportera \`a [HMDP3].
 
\vfill\eject

\titre {0. Notations g\'en\'erales}

Si $A$ est un  anneau, on pose $T= \Spec A$ et on note $\bP^3_A$ ou $\bP^3_T$ l'espace
projectif de dimension
$3$ sur
$T$. On note $R_A$ l'anneau $A[X,Y,Z,T]$. Si
$\cF$ est un faisceau coh\'erent sur
$\bP^3_A$ on note
 $H^i \cF$ le $A$-module $ H^i(\bP^3_A, \cF)$ et on pose $H^i_* \cF=\bigoplus _{n\in
\bZ} H^i \cF (n)$.

Soit $\cN$ un faisceau coh\'erent sur $\bP^3_A$ et
$t$  un point de $T$. On pose $\cN_t= \cN \T_A k(t)$. Si  $u : \cE \fl \cF$ est un
morphisme de faisceaux on
 pose  $u_t=u\T_A k(t)$.

Si $A$ est local, un faisceau
sur $\bP^3_A$  sera dit {\bf dissoci\'e} s'il est de la forme 
$ \cF = \bigoplus_{i=1}^r \cO_{\bP_A}(-n_i) $ o\`u les $ n_i$ sont des entiers.
Si $t$ est le point ferm\'e de $\Spec A$ et si $\cN$ est plat sur $A$ il r\'esulte de [H]
3.10 et 7.9 que
$\cN$ est dissoci\'e si et seulement si $\cN_t$ l'est.

Si $L$ est un $R_A$-module gradu\'e libre  de type fini (resp.
$\cL$ un faisceau dissoci\'e), on
peut l'\'ecrire sous la forme 
$L =
\bigoplus_{n
\in
\bZ} R_A(-n)^{l(n)}$ (resp. $\cL =
\bigoplus_{n
\in
\bZ} \cO_{\bP_A}(-n)^{l(n)}$) o\`u $l$ est une fonction de $\bZ$ dans $\bN$, \`a support
fini, appel\'ee {\bf fonction caract\'eristique} de $L$ (resp. $\cL$). Le plus grand entier
$n$ tel que
$l(n)$ soit non nul est not\'e $\sup L$ (resp. $\sup \cL$).

Si  $q$ est une fonction de $\bZ$ dans $\bZ$  nulle pour
$n \ll 0$, on pose 
$q\sh(n) = \sum_{k \leq n} q(k).$ La donn\'ee de $q$ \'equivaut \`a celle de $q\sh$ en
vertu de la formule 
$q(n) = q\sh (n)- q\sh (n-1)$.

Si $M$ est un $R_A$-module gradu\'e et $n$ un entier on note $M_{\leq n}$ le sous-module
gradu\'e de $M$ engendr\'e par ses \'el\'ements de degr\'e $\leq n$. Si $L$ est libre de
fonction caract\'eristique $l$ on a $\rg L_{\leq n} = l\sh (n)$.

Sur un corps $k$, on appelle courbe un sous-sch\'ema de $\bP^3_k$ de dimension $1$,
sans composante ponctuelle (immerg\'ee ou non), c'est-\`a-dire localement de
Cohen-Macaulay. 

 Une famille de courbes gauches  $\cC$  sur
$T$ est un sous-sch\'ema
ferm\'e de
$\bP^3_T$, plat sur
$T$, et dont les fibres sont des courbes au sens pr\'ec\'edent.

\titre {1. Fonction $q$ relative \`a un sous-module de sections}

\tarte {a) Notations}

Dans ce premier paragraphe on travaille sur un corps $k$ infini. On pose
$\bP^3=\bP^3_k$ et $R=R_k$. Si $x$ est un point de $\bP^3$ et $u : \cE \fl \cF$ un morphisme de faisceaux
coh\'erents on note
$u(x)$ l'application lin\'eaire de $\cE \T_{\cO_\bP} k(x)$ dans $\cF \T_{\cO_\bP} k(x)$
induite par
$u$. Le rang de $u$ est par d\'efinition le rang de $u(\xi)$ o\`u $\xi$ d\'esigne le
point g\'en\'erique de $\bP^3$.
Les notations suivantes seront fix\'ees dans
tout le paragraphe 1 : 

\th {Notations 1.0}. On consid\`ere
 un faisceau coh\'erent sans torsion $\cN$ sur
$\bP^3$, on pose $N = H^0_* \cN$ et on se donne  un sous-$R$-module gradu\'e
$N'$ de 
$N$. Si $n$ est un entier, on  note $\s_{\cN,N',n}$ ou simplement, s'il n'y a pas
d'ambigu\"\i t\'e,
$ \s_n : N'_n \T_k
\cO_{\bP}(-n)
\fl
\cN$ le morphisme canonique.  On note $\cN_{N', \leq n}$ son image qui est aussi le
sous-faisceau de
$\cN$ engendr\'e par
$N'_{\leq n}$.

\th {D\'efinition 1.1}. Avec les notations de 1.0 on d\'efinit la fonction
$q_{\cN,N'}$ (ou simplement $q$ s'il n'y a pas d'ambigu\"\i t\'e)  comme suit : $q\sh
(n)$ est le plus grand entier $m$ tel qu'il existe un faisceau dissoci\'e 
$\cP \subset \cN$, de rang $m$, v\'erifiant  $H^0_* \cP \subset N'$ et $\sup \cP \leq n$,
tel que
$\cN / \cP$ soit sans torsion.

\remas {1.2} \lign
0) On notera que les entiers $n$ qui v\'erifient les conditions de 1.1 forment un intervalle
$] - \infty, m]$. \lign 1) Il est clair que $q\sh$ est une fonction croissante, donc que $q(n)
= q\sh (n)-q\sh (n-1)$ est $\geq 0$. De plus $q\sh$ est nulle pour $n \ll 0$. En effet, comme
$\cN$ est sans torsion, on a $H^0 \cN (n)=0$ pour $n \ll 0$ (cf. [HMDP1] 0.2). Enfin, comme
$q\sh$ est born\'ee par le rang de
$\cN$ elle est constante pour
$n
\gg 0$. On en d\'eduit que
$q$ est
\`a support fini.\lign 2) Par rapport \`a des r\'esultats ant\'erieurs (cf. [MDP1,3]), la
diff\'erence essentielle r\'eside  dans le fait que le faisceau $\cP$ est ``\`a valeurs
dans
$N'$''. Il y a essentiellement deux cas d'application de cette situation :

a) le cas banal, $N=N'$, qui redonne la fonction $q$ ordinaire, voir \S\ {\sl d)} le
calcul explicite de $q$,

b) le cas ``\`a param\`etres'' qui fait l'objet du \S 2 : on part d'un faisceau $\cN$ sur
$\bP^3_A$, o\`u $A$ est un anneau local de corps r\'esiduel $k(t)$,  on consid\`ere le
faisceau $\cN_t= \cN \T_A k(t)$ et on prend, comme sous-module $N'$ de $H^0_* \cN_t$
l'image canonique de $(H^0_* \cN)\T_A k(t)$.

\th {Proposition 1.3}. On reprend les notations de 1.0. Soit $x \in \bP^3$ un point de
codimension
$\leq 1$. On a, pour tout $n \in \bZ$,  $ q\sh
(n) \leq \rg \s_n(x)$.

\dem Par d\'efinition de $q$ il existe une suite exacte $0 \fl \cP \Fl{u} \cN \fl \cE \fl
0$ avec $\sup \cP \leq n$, $\rg \cP = q\sh (n)$, $H^0_* \cP \subset N'$ et $\cE$ sans
torsion. Si
$x$ est un point de codimension $\leq 1$, $\cE$ est plat en $x$ de sorte que la suite
reste exacte en tensorisant par $k(x)$. Comme $\cP$ est en degr\'es $\leq n$ l'image de
$u(x)$ est contenue dans celle de $\s_n (x)$ et on en d\'eduit le r\'esultat.

\tarte {b) Le th\'eor\`eme principal}

Soit $H$  un $k$-espace
vectoriel de
dimension finie. Nous dirons qu'une propri\'et\'e $P$ des \'el\'ements de $H$ est
vraie pour $h$
``g\'en\'eral'' s'il existe un ouvert de Zariski non vide
$U $ du sch\'ema affine associ\'e \`a $H$ tel que $P$ soit vraie pour tout 
$h \in  U$. Puisque $k$ est infini,
 un tel ouvert a des points rationnels.

\th {Lemme 1.4}. On reprend les notations de 1.0. Soit $n \in \bZ$. On pose $\s = \s_n$. Les
assertions suivantes sont \'equivalentes :
\lign 1) Pour $s$ g\'en\'eral dans $N'_n$ on a une suite exacte 
$0 \fl \cO_\bP(-n) \Fl {s} \cN \fl \cE \fl 0$ avec $\cE$ sans torsion. \lign
2) On a les deux conditions suivantes : \lign
\indent \indent a) Pour tout point $x$ de
codimension
$1$ de
$\bP^3$,  on a
$\rg
\s (x)
\geq 1$, \lign
\indent \indent b) on a --- ou bien $\rg \s \geq 2$,\lign 
\indent \indent\indent  \quad \ --- ou bien $\rg \s = 1$, $\Im \s
\simeq
\cO_\bP(-n)$ et $\cN/\Im \s$ sans torsion.

\dem Quitte \`a changer $\cN$ en $\cN (n)$ on peut supposer $n=0$.  Montrons $1 \impl
2$. L'hypoth\`ese implique $q\sh (n) \geq 1$ et la proposition 1.3 montre alors que
$\rg
\s$ et
$\rg
\s (x)$ sont
$\geq 1$.
 Supposons $\rg
\s=1$.  Comme $\Im \s / \Im s$ s'injecte dans $\cE$ il est sans torsion. Comme $\Im
\s$ et $\Im s \simeq  \cO_\bP$ sont tous deux de rang $1$, ils sont \'egaux et on a la
condition requise.

$2 \impl 1$. Si $\rg \s =1$ on a $\Im \s \simeq \cO_\bP$ et on peut prendre pour $s$ 
n'importe quel \'el\'ement non nul de $N'_0$.

Supposons $\rg \s \geq 2$. On consid\`ere le sous-sch\'ema ferm\'e r\'eduit $Z$ de
$\bP(N'_0) \times \bP^3$ adh\'erence de l'ensemble des couples $(s,x)$ tels que
$s(x)=0$. Si on note $d$ la dimension de $\bP(N'_0)$ il suffit de montrer qu'on a $\dim
Z \leq d+1$. En effet, pour $s$ g\'en\'eral dans $N'_0$ la fibre $Z(s)$ sera de
dimension $\leq 1$ ce qui assurera que $\Coker s$ est localement libre en codimension $1$,
donc sans torsion, cf. [MDP1] IV 1.2.

Au point g\'en\'erique $\xi$ de $\bP^3$ on a $\rg \s(\xi) \geq 2$ donc la fibre
$Z(\xi)$ est de dimension $\leq d-2$. En un point $x$ de codimension $1$ on a $\rg
\s(x) \geq 1$ donc 
$\dim Z(x) \leq d-1$. Enfin, en un point de codimension $\geq 2$ les fibres sont de
dimension $\leq d$. En d\'efinitive on a bien $\dim Z \leq d+1$.

\vskip 0.3 cm

\th {Th\'eor\`eme 1.5}. On reprend les notations de 1.0. Soit $\cP$ un faisceau
dissoci\'e de fonction caract\'eristique $p$. Soit $H$ le sous-espace vectoriel 
de  $\Hom_{\cO_{\bP}} (\cP, \cN)$ form\'e des
homomorphismes ``\`a valeurs dans $N'$'' (i.e. tels que l'on ait $u(H^0_* \cP) \subset
N'$).   Les conditions suivantes sont
\'equivalentes : \lign
1) Pour
$u$ g\'en\'eral dans
$H$, on a une suite exacte  $0 \fl \cP \Fl{u} \cN \fl \cE \fl 0$ 
 avec $\cE$ sans torsion. \lign
2) La fonction $p$ v\'erifie : \lign
\indent i) $p\sh (n) \leq q\sh (n)$ pour tout $n$, \lign
\indent ii) s'il existe $n \in \bZ$ avec $p\sh (n) = q\sh (n)= \rg \s_n$, on a un
isomorphisme $\cP_{\leq n} \simeq \cN_{N', \leq n}$.

\dem $1) \impl 2)$. Comme $\cN/ \cP_{\leq n}$ est sans torsion la condition {\sl i)}
r\'esulte de la d\'efinition de la fonction $q$.

Montrons l'assertion {\sl ii)}. Soit $n$ v\'erifiant $p\sh (n) = q\sh
(n) =Ê\rg \s_n$. On a
 le diagramme commutatif de suites exactes
 $$\matrix{
0& \fl & \cP_{\leq n}& \Fl {u}& \cN_{N', \leq n}&
\fl & \cF& \fl& 0\cr
&& \parallel&& \vf && \Vf{j} \cr
0& \fl & \cP_{\leq n} & \fl& \cN&
\fl &\cN/
\cP_{\leq n}
& \fl& 0\cr}$$ 
 Comme $j$ est injective, le
faisceau
$\cF =
\Coker u$ est  sans torsion. Mais, les faisceaux $ \cP_{\leq n}$ et 
 $\cN_{N', \leq n}$ sont tous deux  de rang $p\sh (n) = q\sh (n) =
\rg \s_n$, donc $\cF$ 
est nul  et
$u$ est un isomorphisme. 

$2) \impl 1)$. On proc\`ede par r\'ecurrence sur le rang de $\cP$. Le cas $\rg \cP=0$ est
trivial.

Supposons  $\rg \cP \geq 1$. 
 Soit
$a=
\sup
\cP$. Si on a
 $p\sh (a) = q\sh (a) = \rg \s_a$, la condition {\sl ii)} montre qu'on a un
isomorphisme
$\cP\simeq\cN_{N',\leq a}$. De plus, le quotient $\cN/ \cN_{N', \leq a}$ est sans
torsion. En effet, par d\'efinition de $q \sh (a)$, il existe un sous-faisceau 
dissoci\'e $\cP'$ de $\cN$, de rang $q\sh (a)$, avec $\sup \cP' \leq a$ et $\cN/\cP'$
sans torsion. On a donc une injection 
$j :
\cP' \fl
\cN_{N',\leq a}$. Comme ces faisceaux ont m\^eme rang, $\Coker j$ est de torsion. Mais,
par le lemme du serpent, $\Coker j$ s'injecte dans $\cN / \cP'$ qui est sans torsion,
donc il est nul et $j$ est un isomorphisme.
 Si
$u$  est un
\'el\'ement g\'en\'eral de $H$, $u$ induit un isomorphisme de $\cP$ sur
$\cN_{N',\leq a}$, donc  $\Coker u$ est sans torsion, d'o\`u la conclusion.

\vskip 0.2 cm

Ce cas \'etant maintenant exclu on consid\`ere la fonction $p'$ d\'efinie par $p'(n)=
p(n)$ pour $n<a$ et $p'(a) = p(a)-1$. Cette fonction v\'erifie les m\^emes conditions que
$p$ de sorte que si on pose $\cP' = \bigoplus_{n \in \bZ} \cO_{\bP} (-n)^{p'(n)}$ on
peut appliquer l'hypoth\`ese de r\'ecurrence \`a $\cP'$. Si $u'$ est  un homomorphisme
g\'en\'eral \`a valeurs dans $N'$ on a donc une suite exacte :
$$0 \fl \cP' \Fl {u'} \cN \fl \cE' \fl 0$$
avec $\cE'$ sans torsion. On note $E'$ l'image de $N'$ dans $H^0_* \cE'$ et on
applique le lemme 1.4
\`a
$\cE'$ et \`a l'entier $a$. En un point $x$ de codimension $\leq 1$ de $\bP^3$ le
faisceau
$\cE'$ est localement libre et on a le diagramme commutatif de suites exactes :
$$\matrix {0& \fl & H^0 \cP'(a) \T k(x) & \fl & N'_a\T k(x) & \fl & E'_a \T k(x) & \fl &
0
\cr
&& \vf&&\Vf{\s(x)}&&\Vf{\s'(x)} \cr
0& \fl& \cP' \T k(x) & \fl & \cN \T k(x) & \fl & \cE' \T k(x)& \fl & 0 \cr
} $$
o\`u l'on a pos\'e $\s = \s_{\cN, N',a}$ et $\s' = \s_{\cE', E',a}$.

En un tel point $x$  on a, en vertu de 1.3, $\rg \s(x) \geq q\sh (a)
\geq p
\sh (a)$ et, comme $\cP'$ est de rang $p\sh (a) -1$, on en d\'eduit $\rg \s'(x) \geq 1$.
Si le rang de $\s'$ \'etait \'egal \`a $1$ on aurait $\rg \s = q\sh (a) = p \sh (a)$ or
c'est le cas exclu ci-dessus.

 On a donc $\rg \s' \geq 2$ et, en vertu de 1.4 on obtient, pour $v'$ g\'en\'eral
dans
$E'_a$, une suite exacte 
$$0 \fl \cO_{\bP}(-a) \Fl {v'} \cE' \fl \cE \fl 0$$
avec $\cE$ sans torsion ; la fl\`eche $v'$ se rel\`eve  en $v: \cO_{\bP}(-a) \fl \cN$
\`a valeurs dans $N'$ et si on pose 
$u= v+ u'$, $u$ est g\'en\'eral et on a $\Coker u= \cE$ comme annonc\'e.

\rema {1.6} La condition {\sl 2) ii)} du th\'eor\`eme signifie que le faisceau
$ \cN_{N', \leq n}$ (isomorphe \`a $\cP_{\leq n}$)  est ``obligatoire'' : c'est l'unique
sous-faisceau
$\cL$ de
$\cN$, dissoci\'e et de
 rang
$q\sh (n)$, qui satisfait les conditions de 1.1  :  $\sup \cL \leq n$ et $\cN/ \cL$ sans
torsion.

\th {Proposition 1.7}. On suppose $\cN$ non dissoci\'e et de rang $r$. \lign
1) On a, pour tout $n$, $q \sh (n) \leq r-1$, \lign
2) Si on suppose, de plus, que $N'_n$
engendre
$\cN$ pour
$n
\gg 0$  on a  $q\sh (n) = r -1$ pour $n \gg 0$.

\dem 1) Il est clair qu'on a
$q\sh (n) \leq r-1$. En effet, si $\cP$ est un sous-faisceau dissoci\'e de rang $r$ de
$\cN$ (n\'ecessairement distinct de $\cN$), le quotient a de la torsion.

2) On raisonne par r\'ecurrence sur
$r$. Pour
$r=1$ il n'y a rien
\`a d\'emontrer. Si
$r
\geq 2$ on prend $n \gg 0$ et on applique le lemme 1.4 \`a $\cN$ et $n$. Par
hypoth\`ese on a
$\rg
\s_n(x)
\geq 2$ en tout point, d'o\`u une section $s$ avec un quotient $\cE$ sans torsion. Le
module image $E' \subset H^0_* \cE$ satisfait les m\^emes hypoth\`eses que $N'$. On a
donc, par l'hypoth\`ese de r\'ecurrence, $q\sh_{\cE, E'} (n) = r-2$ pour $n$ grand,
donc un sous-faisceau dissoci\'e $\cP \subset \cE$, de rang $r-2$, \`a valeurs dans
$E'$ et \`a quotient sans torsion. La fl\`eche $\cP \fl \cE$ se rel\`eve \`a $\cN$ et 
la consid\'eration du sous-faisceau $\cO_\bP (-n) \oplus \cP$ de $\cN$ montre qu'on a
$q\sh_{\cN, N'}(n) \geq r-1$.

\tarte {c) Le cas des courbes}

Si on part d'un faisceau $\cN$ r\'eflexif et non dissoci\'e   on obtient le corollaire
suivant   :

\th {Corollaire 1.8}. On suppose que $\cN$  est r\'eflexif et non dissoci\'e    de rang
$r$. Soit $\cP$ un faisceau dissoci\'e
de rang $r-1$ et de fonction caract\'eristique $p$.  Soit $H$ le sous-espace vectoriel de
  $\Hom_{\cO_{\bP}} (\cP, \cN)$ form\'e des
homomorphismes \`a valeurs dans $N'$. Les conditions suivantes sont \'equivalentes : \lign
i) Pour $u$ g\'en\'eral dans  $ H$, on a une suite exacte :
$$0 \fl \cP \Fl {u} \cN \fl \cJ_{C}(h) \fl 0 \leqno {(*)}$$
 o\`u $C$
est une  courbe de $\bP^3$ et $h$ un entier. \lign
ii) La fonction $p$ v\'erifie les deux conditions suivantes : \lign
\indent 1) on a $p\sh (n) \leq q\sh (n)$ pour tout $n \in \bZ$,
\lign
\indent  2) s'il existe $n$ tel que l'on ait $p\sh (n)=q\sh (n)= \rg \s_n$  on a un
isomorphisme $\cP_{\leq n}\simeq \cN_{N',\leq n}$. \lign
Il existe des fonctions $p$ v\'erifiant les conditions {\sl i)} et {\sl ii)} ci-dessus si et
seulement si on a $q\sh (n) = r-1$ pour $n \gg 0$. Dans ce cas, la fonction $q$ elle-m\^eme
convient et la valeur minimale du d\'ecalage
$h$ correspond au cas $p=q$ et est
\'egale
\`a
$\sum_{n\in
\bZ} nq(n) +
\deg
\cN$.

\dem Le sens $i) \impl ii)$ vient de 1.5. R\'eciproquement,  on a, en vertu de 
1.5, un homomorphisme injectif
$u :
\cP
\fl
\cN$ dont le conoyau
$\cE$ est un faisceau sans torsion de rang $1$ donc l'id\'eal tordu $\cJ_Z (h)$ d'un
sous-sch\'ema $Z$ de $\bP^3$ de dimension $\leq 1$. Comme $\cN$ n'est pas dissoci\'e,
$Z$ n'est pas vide. Par ailleurs, comme
$\cN$ (resp. $\cP$) est r\'eflexif  on a ${\rm prof}\, \cJ_{Z,x} \geq 2$ en tout point
ferm\'e de $\bP^3$,  donc ${\rm prof}\, \cO_{Z,x} \geq 1$ et, {\it a fortiori} $\dim 
\cO_{Z,x} \geq 1$, de sorte que $Z$ est une courbe.

S'il existe une fonction $p$ v\'erifiant les conditions ci-dessus, on a $p\sh (n)=
r-1$ pour $n \gg 0$ donc aussi $q\sh (n) = r-1$ pour $n \gg 0$ par {\sl ii)} et 1.7.
R\'eciproquement, si on a cette condition, il existe,  par d\'efinition de $q$, un
sous-faisceau
$\cP$ de $\cN$, dissoci\'e de fonction caract\'eristique $q$  dont  le quotient
$\cN/
\cP$ est sans torsion de rang $1$, donc, par le m\^eme argument que ci-dessus, c'est un
id\'eal tordu,  de sorte que
$q$ convient. Enfin,  la condition
$p\sh
\leq q\sh$, jointe
\`a la formule 
$\sum_{n \in \bZ} nq(n) = - \sum_{n \in \bZ} q \sh (n)$, montre que le d\'ecalage
est minimal dans le cas de la fonction $q$.

\remas {1.9} \lign
1) L'existence d'un faisceau dissoci\'e  $\cP$ v\'erifiant les conditions de 1.8 d\'epend
essentiellement du sous-module $N'$.  Pr\'ecis\'ement, on a vu ci-dessus que cette existence \'equivaut
\`a la condition
$q\sh_{\cN,N'} (n)= r-1$ pour
$n \gg 0$.  D'apr\`es 1.7  cette condition est satisfaite si $N'_n$ engendre $\cN$ pour $n
\gg 0$. Dans certains cas il n'y a pas de faisceau $\cP$ v\'erifiant les conditions de
1.8 (par exemple si
$N'$ est nul et
$\cN$ de rang
$\geq 2$).
\lign 2) Si le faisceau $\cN$ est dissoci\'e le corollaire est en d\'efaut lorsque $u$ (\`a
valeurs dans $N'$) identifie
$\cP$ \`a  un facteur direct de $\cN$  car le quotient est alors un faisceau $\cO_\bP (n)$. Dans le
cas
$N'=N$ la fonction $q$ est la fonction caract\'eristique de $\cN$ et on v\'erifie facilement
que le seul cas
\`a
\'ecarter est celui o\`u l'on a $p(n) \leq q(n)$ pour tout
$n$.

\tarte {d) Calcul de la fonction $q$}

\th {D\'efinition 1.10}. Avec les notations de 1.0 on pose $\a_n =
\a_n (\cN, N') = \rg \s_n$ et
$$\beta_n= \beta_n (\cN, N') = {\rm Min}\, \{ \rg \s_n (x) \mud x \in \bP^3, x \ \hbox
{de codimension 1}\}.
$$ 
\rema {1.11} On a $\beta_n \leq \a_n \leq \rg \cN$. En vertu de 1.3 on a aussi $q\sh  (n)
\leq
\beta_n \leq 
\a_n$.

\vskip 0.3 cm
La proposition suivante \'etudie le cas d'\'egalit\'e  $q\sh (n) = \a_n$ et 
 le ph\'enom\`ene d\'ej\`a
rencontr\'e  en   1.5 (l'existen\-ce d'une partie ``obligatoire'' pour $\cN$) :

\th {Proposition 1.12}. Soit $a \in \bZ$. On suppose qu'on a $q\sh (a) = \a_a$,
c'est-\`a-dire qu'il existe un sous-faisceau dissoci\'e 
$\cL
\subset
\cN$ v\'erifiant $H^0_* \cL \subset N'$, $\sup \cL \leq  a$, $\rg \cL = \a_a=\rg \s_a$ et
$\cN/
\cL$ sans torsion. Alors,  on a n\'ecessairement  $\cL= \cN_{N',\leq a}$ et, en posant $L =
H^0_* \cL$, $L= N'_{\leq a}$. Si
$l$ est la fonction ca\-rac\-t\'eristique de $\cL$ on a $l\sh (n) = \a_n=
\beta_n = q\sh (n)$ pour $n \leq a$ et $l\sh (n) = \a_a$ pour $n>a$.

\dem On pose $\cK= \cN/\cL$ et le lemme du serpent
montre que dans  le diagramme  :
$$\matrix {
 0 & \fl & \cL& \fl& \cN& \fl& \cK& \fl& 0 \cr
&&\Vf{j}&&\parallel&&\Vf{p} \cr
0 & \fl & \cN_{N', \leq a}& \fl& \cN& \fl& \cK_{N',a}& \fl& 0 \cr
}$$
le noyau $\cF$ de $p$ est isomorphe au conoyau de $j$. Comme $\cK$ est sans torsion il en
est de m\^eme de
$\cF$, ce qui, puisque
$\cL$ et
$\cN_{N',
\leq a}$ sont tous deux de rang
$\a_a$, montre que  
$j$ est   un
isomorphisme. Comme $L= H^0_* \cL$ est inclus dans
$N'$, donc aussi dans $N'_{\leq a}$, cela montre aussi l'\'egalit\'e $L= N'_{\leq a}$. 

Comme $L$ est libre, on en d\'eduit, pour $n \leq a$, $  L_{\leq n}= N'_{\leq n}$, d'o\`u
$l\sh (n) =
\rg L_{\leq n} = \rg N'_{\leq n} = \a_n= \beta_n$.  Il est clair que $\cN/ \cL_{\leq n}$ est
sans torsion et on a donc aussi $\a_n= q\sh (n)$. Enfin, comme $\sup L \leq a$, on a $l\sh
(n) = l\sh (a) =
\a_a$ pour
$n>a$.

\vskip 0.2 cm

Comme $\cN$ est sans torsion tout entier $n \ll 0$ v\'erifie la condition $q\sh (n) =
\a_n=0$ et on peut donc poser la d\'efinition suivante :

\th {D\'efinition 1.13}. Avec les notations de 1.0, on d\'efinit l'entier $b_0= b_0
(\cN, N')$
\footnote {$(^1)$}{Par rapport aux notations usuelles de [MDP] 1,3 on a pos\'e $b_0=
a_0-1$.} comme le plus grand entier
$n$ qui v\'erifie $q\sh (n)  = \rg \s_n = \a_n$. Si tous les entiers $n$ conviennent
on pose $b_0= +
\infty$.

\vskip 0.3 cm 

La proposition suivante, qui est imm\'ediate \`a partir de 1.12, sera utile pour le calcul
explicite de
$b_0$ :

\th {Proposition 1.14}.  Soit  $a \in \bZ$. Les conditions suivantes sont
\'equivalentes : \lign
i) on a  $a \leq b_0$, \lign
ii) on a $q\sh (a) = \a_a$, \lign
iii) il existe un sous-faisceau dissoci\'e 
$\cL
\subset
\cN$ v\'erifiant $H^0_* \cL \subset N'$, $\sup \cL \leq  a$, $\rg \cL = \rg \s_a$ et
$\cN/
\cL$ sans torsion, \lign
iv)
 le sous-module $N'_{\leq a}$ de $N'$  est libre de rang $\a_a$ et le quotient
$N/N'_{\leq a}$ est sans torsion.

\th {Th\'eor\`eme 1.15}. La fonction $q$ se calcule comme suit :
\lign 1) pour $n \leq b_0$ on a $q\sh (n) = \a_n = \beta_n$, \lign
2) pour $n > b_0$ on a $q\sh (n) = \inf \, (\a_n -1, \beta_n)$.

\dem L'assertion 1) r\'esulte de 1.12. Supposons $n > b_0$ (donc $\cN \neq 0$). On a  $q\sh
(n) <
\a_n$ par d\'efinition de $b_0$ et 
$q\sh (n)
\leq
\beta_n$ par 1.11, donc $q\sh (n) \leq \inf \, (\a_n -1, \beta_n)$. R\'eciproquement,
posons $d_n=d_n(\cN,N') = \inf \, (\a_n -1, \beta_n)$ et montrons, par r\'ecurrence sur
le rang de $\cN$, qu'on a $q\sh (n) \geq d_n$, c'est-\`a-dire qu'on a un sous-faisceau
dissoci\'e
$\cP$ de $\cN$, en degr\'es $\leq n$,  de rang $d_n$, \`a quotient sans torsion. Si
$d_n$ est nul il n'y a rien
\`a montrer. Supposons $d_n \geq 1$. En vertu de 1.4 il existe une suite exacte
$0 \fl \cO_\bP (-n) \Fl{s} \cN \fl \cE \fl 0$ avec $s$ \`a valeurs dans $N'$ et  $\cE$
sans torsion. Soit
$E'$ l'image de $N'$ dans $H^0_* \cE$. On a  le lemme suivant :

\th {Lemme 1.16}. On a $d_n(\cE, E') \geq d_n (\cN,N')-1$.

\dem (de 1.16) Soit $x$ un point de codimension $\leq 1$ de $\bP^3$. On a le diagramme
de suites exactes suivant : 
$$\matrix {0& \fl &  k(x) & \fl & N'_n\T_k k(x) & \fl & E'_n \T_k k(x) & \fl &
0
\cr
&& \vf&&\Vf{\s_n(x)}&&\Vf{\s'_n(x)} \cr
0& \fl& \cO_\bP(-n) \T_{\cO_\bP} k(x) & \fl & \cN \T_{\cO_\bP} k(x) & \fl & \cE
\T_{\cO_\bP} k(x)& \fl & 0
\cr } $$
dont on d\'eduit l'in\'egalit\'e $\rg \s'_n(x) \geq \rg \s_n (x)-1$, d'o\`u $ \a_n
(\cE, E')  \geq \a_n (\cN,N')
-1$ et $ \beta_n (\cE, E') \geq \beta_n (\cN,N') -1$
 et la conclusion du lemme.

\vskip 0.3 cm

On peut alors finir de prouver 1.15. En vertu de l'hypoth\`ese de r\'ecurrence on a,
dans tous les cas, 
$q\sh_{\cE,E'} (n) \geq d_n (\cE,E') \geq d_n (\cN,N')-1$. Par
d\'efinition de la fonction
$q$ cela signifie qu'il existe un sous-faisceau dissoci\'e $\cP'$ de $\cE$, de rang
$\geq d_n (\cN,N')-1$, en degr\'es $\leq n$, avec $\cE/\cP'$ sans torsion. La fl\`eche
$\cP' \fl
\cE$ se rel\`eve
\`a $\cN$ et
$\cO_\bP(-n) \oplus \cP'$ fournit le faisceau $\cP$ cherch\'e.

\titre {2. Construction de familles de courbes}

Dans tout ce paragraphe on travaille sur  un anneau noeth\'erien local $A$ d'id\'eal
maximal $m_A$. On pose
$T =
\Spec A$ et on d\'esigne par
$t$ le point ferm\'e de $T$ et par $k(t)$ le corps r\'esiduel. On suppose
$k(t)$ infini, cf. 2.11 pour une variante lorsque $k(t)$ est fini.

\vskip 0.3 cm

\tarte {a) R\'esultats pr\'eliminaires }

\th {D\'efinition 2.1}. On dit que le faisceau
$\cN$ est sans torsion dans la fibre sp\'eciale si  le faisceau 
$\cN_t$ est sans torsion. 
\vskip 0.3 cm

\th {Proposition 2.2}. 
On suppose qu'on a une suite exacte $\cF \Fl {u} \cG \fl \cH \fl 0$ de faisceaux
coh\'erents sur $\bP^3_A$ avec
$\cG$ plat sur
$A$ \footnote {$(^2)$} {Dans ce paragraphe, la platitude est toujours sous-entendue sur
$A$.} et
$\cF$ dissoci\'e. Alors, les conditions suivantes sont
\'equivalentes :
\lign a) $u$ est injectif et $\cH$ est
plat sur $A$, \lign
b)  $u_t: \cF_t \fl \cG_t$ est injectif au point ferm\'e $t$ de $T$.

\dem

 Cela r\'esulte du crit\`ere local de platitude, cf. [EGA] $0_{III}$ 10.2.4 et de [EGA] IV,
11, Note de la page 118.

\vskip 0.3 cm

Nous faisons maintenant le lien avec la situation du paragraphe 1.

\vskip 0.2 cm

Soit $\cN$ un faisceau coh\'erent sur $\bP^3_A$ et soit $\rho : H^0_* (\bP^3_A,\cN) \T_A k(t)
\fl H^0_* (\bP^3_{k(t)},\cN_t)$ l'homomorphisme canonique. On note $N'$ l'image de $\rho$
dans $N = H^0_* \cN_t$. En g\'en\'eral $N'$ est diff\'erent de $N$ (dire que $N'=N$
signifie que $H^0_* \cN$ commute au changement de base, cf. [AG] III 12.11). Dans toute la
suite nous utilisons les r\'esultats du paragraphe 1 appliqu\'es au faisceau $\cN_t$ et
au sous-module $N'$. On notera que, si $\cN$ est plat sur $A$ et $n \gg 0$, $H^0 \cN (n)$
commute au changement de base, donc que $N'_n$ engendre
$\cN_t$, cf. 1.7.

\th {Corollaire 2.3}. Soit $\cN$ un faisceau coh\'erent sur $\bP^3_A$, plat et sans
torsion dans la fibre sp\'eciale,  soit $N' =
\Im
\rho$ comme ci-dessus et soit $\cP$ un faisceau dissoci\'e sur $\bP^3_A$ de fonction
caract\'eristique
$p$.  Les propri\'et\'es suivantes sont \'equivalentes : \lign 1) Il existe une suite exacte
de faisceaux sur $\bP^3_A$
$$0 \fl \cP \Fl {u} \cN \fl \cE \fl 0$$ avec  $\cE$ plat et sans
torsion dans la fibre sp\'eciale. \lign
2) Il existe une suite exacte de faisceaux sur $\bP^3_{k(t)}$
$$0 \fl \cP_t \Fl{u_0} \cN_t \fl \cE_0 \fl 0$$ avec $u_0$
 \`a valeurs dans $N'$ et $\cE_0$  sans torsion.

\dem Le sens $1) \impl 2)$ est imm\'ediat en tensorisant par $k(t)$ puisque $\cE$ est
plat sur $A$ (comme $u_t$  provient de $u$ il est bien \`a valeurs dans $N'$). Dans
l'autre sens, comme
$P_t = H^0_*
\cP_t$ est libre sur
$R$, la fl\`eche
$P_t
\fl N'$ se rel\`eve en
$P_t
\fl H^0_* \cN \T_A k(t)$ puis en $P \fl H^0_* \cN$  et  induit $u : \cP \fl \cN$ avec
$u_t=u_0$. On conclut avec 2.2.

\tarte {b) Construction de familles de courbes}

Nous d\'efinissons  les invariants associ\'es \`a un faisceau $\cN$ sur $\bP^3_A$ ; ce
sont ceux du faisceau $\cN_t$ relatifs au sous-module $N'$ :

\th {D\'efinition 2.4}. Soit $\cN$ un faisceau coh\'erent  sur $\bP^3_A$, plat et sans
torsion dans la fibre sp\'eciale. On consid\`ere   l'homomorphisme canonique
$\rho : H^0_* (\bP^3_A,\cN)
\T_A k(t)
\fl H^0_* (\bP^3_{k(t)},\cN_t)$  et on note $N'$ son image. \lign
On pose alors $\cN_{t, \leq n} = (\cN_t)_{N', \leq n}$,
$\a_n (\cN) = \a_n (\cN_t, N')$, $\beta_n (\cN) = \beta_n (\cN_t, N')$, $b_0 (\cN) =
b_0 (\cN_t, N')$,
$q_\cN = q_{\cN_t, N'}$. S'il n'y a pas d'ambigu\"\i t\'e sur le faisceau $\cN$ ces
invariants seront not\'es simplement
$\a_n$, $\beta_n$, $b_0$, $q$.

Soient $H$  un $A$-module de type fini, $\ov H =H \T_A k(t)$ le $k(t)$-espace
vectoriel de
dimension finie obtenu par r\'eduction modulo $m_A$ et $\pi : H \fl \ov H$
la projection
canonique. Nous dirons qu'une propri\'et\'e $P$ des \'el\'ements de $H$ est
vraie pour $h$
``g\'en\'eral''dans $H$ s'il existe un ouvert de Zariski non vide
$\ov U $ du sch\'ema affine associ\'e \`a $\ov H$ tel que $P$ soit vraie pour tout 
$h \in \pi^{-1} (\ov U)$. Puisque $k(t)$ est infini,
 un tel ouvert a des points rationnels, et son image r\'eciproque $\pi^{-1} (\ov U)$,
consid\'er\'ee comme sous-ensemble du $A$-module $H$, n'est pas vide.

Avec ces notations et compte tenu de la d\'efinition de $b_0$, le corollaire 2.3,
joint
\`a 1.5, donne aussit\^ot le th\'eor\`eme suivant :

\th {Th\'eor\`eme 2.5}.  On suppose   que le
faisceau $\cN$ est  plat et sans
torsion dans la fibre sp\'eciale. Soit $\cP$ un faisceau
dissoci\'e  de fonction caract\'eristique $p$. 
Les conditions suivantes sont \'equivalentes  : \lign
i) Pour $u$ g\'en\'eral dans $H=\Hom_{\cO_{\bP}} (\cP, \cN)$, on a une suite exacte 
$$0 \fl \cP \Fl {u} \cN \fl \cE \fl 0$$
o\`u $\cE$ est plat et sans
torsion dans la fibre sp\'eciale. \lign
ii) La fonction $p$ v\'erifie les conditions 1) et 2) ci-dessous : \lign \indent 1) on a
$p\sh (n)
\leq q\sh (n)$ pour tout
$n
\in
\bZ$,
\lign
\indent  2) s'il existe $n \leq b_0$ tel que l'on ait $p\sh (n)=q\sh (n)$ on a un
isomorphisme $(\cP_{\leq n})_t \simeq \cN_{t,\leq n}$. \lign
La fonction $q$ v\'erifie les conditions {\sl i)} et {\sl ii)} ci-dessus.

Dans le cas des courbes on a le r\'esultat suivant :

\th {Th\'eor\`eme 2.6}.  Soit $\cN$   un faisceau sur $\bP^3_A$,
localement libre de rang $r$ et non dissoci\'e.  Soit
$\cP =
\bigoplus_{n
\in \bZ} \cO_{\bP_A} (-n)^{p(n)}$ un faisceau dissoci\'e  de rang $r-1$.  Alors,
les conditions suivantes sont \'equivalentes : \lign i) Pour $u$ g\'en\'eral
dans 
$ H=\Hom_{\cO_{\bP}} (\cP,
\cN)$, on a une suite exacte :
$$0 \fl \cP \Fl {u} \cN \fl \cJ_{\cC}(h) \fl 0 $$
 o\`u $\cC$
est une famille de courbes de $\bP^3_A$, plate  sur $A$. \lign
ii) 
La fonction $p$ v\'erifie les conditions suivantes : \lign
\indent 1) on a $p\sh (n) \leq q\sh (n)$ pour tout $n \in \bZ$, \lign
\indent 2) s'il existe $n\leq b_0$ tel que l'on ait $p\sh (n)=q\sh (n)$  on a un
isomorphisme $(\cP_{\leq n})_t \simeq \cN_{t,\leq n}$. \lign
La fonction $q$ v\'erifie les conditions {\sl i)} et {\sl ii)} ci-dessus. La valeur minimale
du d\'ecalage
$h$ est
\'egale
\`a
$\sum_{n\in
\bZ} nq(n) +
\deg
\cN$. La famille de courbes correspondant \`a ce d\'ecalage minimum est donn\'ee par une
r\'esolution comme ci-dessus, avec $p=q$.

\dem Il reste \`a prouver $ii) \impl i)$. En vertu de 2.5 on a un homomorphisme injectif $u
:
\cP
\fl
\cN$ dont le conoyau
$\cE$ est un faisceau de rang $1$, plat et sans
torsion dans la fibre sp\'eciale. Par ailleurs, comme
$\cN$ (resp. $\cP$) est localement libre de rang $r$ (resp. $r-1$) on a un isomorphisme
$\cN(-d)
\simeq \wedge^{r-1} \cN^{\vee}$ (resp. $\cO_{\bP^3_A} (-e) \simeq \wedge^{r-1} \cP^\vee$)
o\`u $d$ (resp. $e$) est le degr\'e de $\cN$ (resp. $\cP$). On en d\'eduit un complexe 
$$ 0 \fl \cP \Fl {u} \cN\ \Fl {\wedge^{r-1} u^\vee (d)}\ \cO_{\bP^3_A} (h)$$
avec $h= d-e$.
Le conoyau de $\wedge^{r-1} u^\vee (d)$ est de la forme $\cO_{\cC}(h)$  o\`u $\cC$ est un
sous-sch\'ema ferm\'e de $\bP^3_A$. Comme $\cE$ est sans torsion dans la fibre sp\'eciale il
r\'esulte du th\'eor\`eme de Hilbert-Burch, cf. par exemple [MDP1] II, 1.1 que le complexe
ci-dessus est une suite exacte dans la fibre sp\'eciale, donc une suite exacte sur $\bP^3_A$
en vertu de 2.2 appliqu\'e \`a la suite $ \cE \Fl {u} \cO_{\bP^3_A} (h) \fl \cO_\cC (h) \fl
0$.  De plus, comme
$\cN$ n'est pas dissoci\'e, 
$\cC$ est non vide donc est  une famille de courbes (cf. [MDP3] III 1.7) plate sur
$A$ en vertu du crit\`ere local de platitude.

Comme $N'_n$ engendre $\cN_t$ pour $n \gg 0$ on a $q\sh (n) = r-1$ pour $n \gg 0$ (cf. 
1.9.1), de sorte que $q$ v\'erifie les conditions de 2.6.  Enfin,  la condition $p\sh \leq
q\sh$, jointe
\`a la formule 
$\sum_{n \in \bZ} nq(n) = - \sum_{n \in \bZ} q \sh (n)$, montre que le d\'ecalage
est minimal dans le cas de la fonction $q$.

\th {Corollaire 2.7}.  Soit $\cN$ un faisceau localement libre de rang $r$ sur
$\bP^3_A$. Il existe une famille $\cC$ de courbes  lisses et connexes
  de
$\bP^3_A$, plate  sur $A$, et un entier $h$ tels que $\cC$ admette une r\'esolution de
type {\sl N} de la forme :
$$ 0 \fl \cP \fl \cN \fl \cJ_\cC (h) \fl 0.$$

\dem On choisit un entier $n$ assez grand pour que $H^0 \cN (n)$ commute au changement de
base et que $\cN (n-1)$ soit engendr\'e par ses sections globales. On pose
$\cP = \cO_\bP(-n-1)^{r-1}$.  On montre comme ci-dessus qu'un homomorphisme g\'en\'eral $u :
\cP \fl \cN$ a pour conoyau l'id\'eal tordu d'une famille plate de courbes $\cC$ (comme les
degr\'es dans $\cP$ ont \'et\'e choisis assez grands il n'est pas n\'ecessaire de supposer
$\cN$ non dissoci\'e ici). De plus, les hypoth\`eses montrent que la courbe sp\'eciale est
lisse connexe (en toute caract\'eristique), cf. [K] Th. 3.3, donc aussi les autres.

\tarte {c) Application \`a la biliaison}

Nous d\'ecrivons maintenant  les familles de courbes qui sont  minimales dans les classes
de biliaison. 

 Rappelons d'abord, cf. [HMDP1] que les classes de biliaison 
correspondent, via les r\'esolutions de type \sN, aux classes de faisceaux localement
libres pour la relation  de pseudo-isomorphisme et qu'une telle classe
contient un  faisceau localement libre $\cN_0$ extraverti (c'est-\`a-dire v\'erifiant 
$\Ext^1_{R_A} (H^0_* \cN_0, R_A)=H^1_* \cN_0^\vee =0$), minimal (i.e. sans
facteur direct dissoci\'e), unique
\`a isomorphisme pr\`es, cf. [HMDP1] 2.14.

Nous aurons besoin aussi de la notion de d\'eformation \`a cohomologie uniforme :

\th {D\'efinition 2.8}. Posons  $S=\Spec A[\la]$ et soit $\Om$ un ouvert de $S$. On note
$A_\Om$  l'ensemble des $a
\in A$ tels que la section $\la= a$ soit contenue dans $\Om$. Soit $\cC \fl \Om$ une
famille plate de courbes. Pour $a \in A_\Om$ on note
$\cC_a$ la famille de courbes (param\'etr\'ee par
$A$) correspondante. Si les points $a$ et $b$ sont dans $A_\Om$  on dit  que $\cC$ est une
d\'eformation joignant
$\cC_a$ et $\cC_b$. On dit que cette d\'eformation est
{\bf \`a cohomologie uniforme } si, pour tout $s \in \Spec A$, tout $i= 0,1,2$ et tout $n \in
\bZ$ la fonction $a \fl h^i(\cJ_{\cC_a}\T_A k(s)) (n)$ est constante sur $A_\Om$ (mais
attention, cette dimension peut varier avec
$s$).

Avec ces notations on a le r\'esultat suivant :

\th {Corollaire 2.9}.  Soit $\cN_0$ un faisceau localement libre
extraverti minimal. Posons   $q=
q_{\cN_0}$ et $h_0= \sum_{n \in \bZ} nq(n) + \deg \cN_0$. \lign
1)  Il existe une
famille de courbes  $\cC_0$ et une r\'esolution 
$ 0 \fl \cP_0 \Fl{v} \cN_0 \fl \cJ_{\cC_0} (h_0) \fl 0$ avec $\cP_0$ dissoci\'e. \lign
2) Si $\cC_1$ est une famille
de courbes  de la classe de biliaison de $\cC_0$, elle admet une r\'esolution 
$ 0 \fl \cP \fl \cN_0 \oplus \cL \fl \cJ_{\cC_1} (h) \fl 0$ avec $\cL$ dissoci\'e et on  on a
$h
\geq h_0$. Si $d$ et $g$ (resp. $d_0$ et $g_0$) sont respectivement le degr\'e et le
genre de $\cC_1$ (resp. $\cC_0$) on a $d \geq d_0$ et $g \geq g_0$.\lign
3) R\'eciproquement, pour tout $h \geq h_0$, il existe une famille de courbes $\cC_1$
avec une r\'esolution comme ci-dessus.
\lign
4) Si  on a
$h=h_0$,
les familles $\cC_0$ et $\cC_1$ sont jointes par une d\'eformation  \`a cohomologie
uniforme. On dit que
$\cC_0$ est une {\bf famille minimale} de courbes de la classe de biliaison associ\'ee
\`a $\cN_0$.

\dem L'existence de $\cC_0$  et de la suite exacte vient de 2.6. Le faisceau $\cP_0$ est
un faisceau dissoci\'e de fonction caract\'eristique $q$.

Soit
$\cC_1$  une  famille de courbes de la classe de biliaison de $\cC_0$. D'apr\`es [HMDP1] 2.22
il existe une r\'esolution de type \sN\ extravertie de $\cC_1$, avec un faisceau $\cN$.
D'apr\`es [HMDP1] 3.2, ce faisceau est pseudo-isomorphe \`a $\cN_0 (h')$ pour un d\'ecalage
$h'$ convenable. Il r\'esulte alors de [HMDP1] 2.14 que l'on a (\`a d\'ecalage pr\`es)
$\cN = \cN_0
\oplus \cL$ avec $\cL$ dissoci\'e de fonction caract\'eristique $l$. On v\'erifie
aussit\^ot la formule
$q_\cN\sh(n) = q_{\cN_0}\sh (n) + l\sh (n)$ (cf. [MDP1] IV 2.9) et on en d\'eduit $h \geq
h_0$. Les assertions sur degr\'e et genre r\'esultent de [MDP1] IV 5.2 et  5.3
appliqu\'es au point $t$.

Le point 3) s'obtient en faisant des biliaisons \'el\'ementaires
triviales, cf. [HMDP1] 1.6, \`a partir de $\cC_0$.

Supposons maintenant $h=h_0$. On a  une suite exacte $0 \fl \cP \Fl {u} \cN_0 \oplus \cL \fl
\cJ_{\cC_1}(h_0) \fl 0$ et, comme $h_0$ est le minimum possible, le faisceau dissoci\'e
$\cP$ a pour fonction caract\'eristique $q_\cN$ (cf. 2.6), donc s'\'ecrit $\cP = \cP_0
\oplus
\cL$ avec $\cP_0$  comme ci-dessus.  

Posons $u= ^t\!(u_1,u_2)$ o\`u $u_2$ va de $\cP_0 \oplus \cL$ dans $\cL$. Nous allons montrer
que $u_2$ admet une section.  Par r\'ecurrence sur le rang de
$\cL$ on se ram\`ene au cas $\cL = \cO_{\bP^3_A} (-a)$ et il faut montrer 
que la matrice $u_2$ a un coefficient qui est une constante inversible. On a la
suite exacte
$$0 \fl (\cP_0)_{\leq a} \oplus \cO_{\bP^3_A} (-a) \Fl{w} \cN_0 \oplus \cO_{\bP^3_A} (-a) \fl
\cE \fl 0$$
o\`u $\cE$, qui est extension de $\cJ_{\cC_1}(h_0) $ par le faisceau dissoci\'e
$\cP_{>a}$, est encore plat et sans
torsion dans la fibre sp\'eciale. Posons $w= ^t\!(w_1,w_2)$ et soit $\cF = \Coker w_1$.
Si
$w_2$ n'a pas de coefficient constant inversible elle est nulle au point ferm\'e $t$ 
et on en d\'eduit que $(w_1)_t$ est injective et que $\cF_t$ est facteur direct de
$\cE_t$, donc sans torsion. Mais,  cela implique $q_{\cN_0}\sh (a) +1 \leq
q_{\cN_0}\sh (a)$, ce qui est absurde.

\vskip 0.2 cm

 L'existence de la section de $u_2$ permet de ``simplifier'' les facteurs $\cL$
c'est-\`a-dire d'obtenir une suite exacte
 $0 \fl \cP_0 \Fl {u} \cN_0 \fl \cJ_{\cC_1} (h_0) \fl 0$.

On consid\`ere alors  
l'homomorphisme $w_\la= \la u + (1-\la) v : \cP_0 \T_A A[\la] \fl
\cN_0\T_A A[\la]$ et  son image $ \ov w_\la = \la \ov u + (1-\la) \ov v$ au point ferm\'e. En
vertu de 2.6 et de la d\'efinition de ``g\'en\'eral'', il existe 
un ouvert affine
$\Omega=
\Spec
\Lambda$ de la droite affine 
$\Spec A[\la]$, contenant les points $0$ et $1$ de la fibre ferm\'ee, donc aussi
les sections $\la=0$ et $\la=1$,  tel que, sur cet ouvert, 
$w_\la$ d\'efinit une famille plate de courbes sur $A$. On consid\`ere alors la
suite exacte de faisceaux sur $\bP^3_\Omega$ :
$$ \cP_{0}\T_A \Lambda \Fl {w_\la} \cN_{0}\T_A \Lambda \fl \cE_\Lambda \fl 0.$$

En vertu de 2.2 appliqu\'e aux localis\'es de $A[\la]$, on voit que $w_\la$ est injectif et
que son conoyau est plat, donc d\'efinit une famille de courbes  $\cC$  param\'etr\'ee par
$
\Omega$, dont les sections correspondant \`a
$\la=0$ et
$\la=1$ sont respectivement $\cC_0$ et $\cC_1$, comme annonc\'e. Comme les familles $\cC_a$
pour $a \in A_\Om$ ont toutes pour r\'esolution $0 \fl \cP_0 \fl \cN_0 \fl \cJ_{\cC_a} (h) \fl
0$, on voit que la d\'eformation est \`a cohomologie uniforme.

 \th {Corollaire 2.10}. (Propri\'et\'e de Lazarsfeld-Rao)  Soit
$\cC_0$ une famille minimale de courbes et $\cC$  une famille de courbes de la classe de
biliaison de
$\cC_0$. Alors, il existe un entier $m \geq 0$ et une suite de courbes $\cC_0, \cC_1,
\cdots, \cC_m$ telle que $\cC_{i+1}$ s'obtienne \`a partir de $\cC_i$ par une
biliaison \'el\'ementaire  (cf. [HMDP1] ) et $\cC$ \`a partir de $\cC_m$ par
une d\'eformation  \`a cohomologie uniforme.

\dem La d\'emonstration originelle de Lazarsfeld et Rao (cf. [MDP1] IV 5.1) s'applique
presque sans changement. Elle repose sur le corollaire 2.9, la seule diff\'erence
\'etant qu'on doit trouver une famille plate de surfaces (au lieu d'une seule surface)
pour effectuer une biliaison. Il suffit pour cela de faire le raisonnement habituel au
point ferm\'e.

\rema {2.11}   Si le corps $k(t)$ est fini, les assertions d'existence de 2.5,
2.6, 2.9 et 2.10 restent vraies \`a condition de remplacer \'eventuellement $A$ par un anneau
local
$B$ fini et \'etale sur $A$. 

\titre {3. Exemples}

\tarte {a) Un algorithme de calcul}

Soit $A$  un anneau local noeth\'erien, $T= \Spec A$ et $t$ le point ferm\'e de $T$. 

Consid\'erons la situation suivante : on se donne deux $R_A$-modules libres gradu\'es
$L_1$ et $L_2$ et un homomorphisme de degr\'e $0$, $s : L_2 \fl L_1$, donn\'e par une
matrice \`a coefficients homog\`enes dans $R_A$ que l'on note encore $s$. Soit
$\wi s :
\cL_2 \fl \cL_1$ la fl\`eche de faisceaux associ\'ee \`a $s$. On pose $\cN =
\Im \wi s$ et $\cK = \Coker
\wi s$ ; on a donc la suite
exacte $(*)$ : $ 0 \fl \cN \fl \cL_1 \fl \cK \fl 0$ et on fait les deux hypoth\`eses
suivantes : 

1)  $\cK$ est localement libre. Cette
condition se traduit  en termes matriciels : si
$r$ est le rang de
$s_t: L_{2,t} \fl L_{1,t}$  (c'est-\`a-dire la dimension d'un plus grand mineur non nul
de la matrice
$s_t$) $\cK$ est localement libre si et seulement si les $r$-mineurs de $s_t$ 
d\'efinissent le vide dans $\bP^3_{k(t)}$.

2) la fl\`eche $L_2 \fl N= H^0_* \cN$ induite par $s$ est surjective.

On notera que tout faisceau localement libre $\cN$ sur $\bP^3_A$ peut s'obtenir de cette
mani\`ere~: pour obtenir $L_1$
 il suffit  de consid\'erer le dual $\cN^\vee$ de $\cN$ et de prendre une
 pr\'esentation libre
$L_1^\vee \fl H^0_* \cN^\vee \fl 0$ et pour obtenir $L_2$ de prendre une pr\'esentation
libre
$L_2 \fl  N= H^0_* \cN \fl 0$. 

Nous montrerons, dans [HMDP3], comment obtenir
syst\'ematiquement des faisceaux
$\cN$ dits triadiques, munis de suites exactes 
$$0 \fl \cN \fl \cL_1 \fl \cL_0 \fl \cL_{-1} \fl 0,$$
avec les $\cL_i$ dissoci\'es qui conduisent naturellement \`a de telles \'ecritures.

\vskip 0.3 cm

On se propose de calculer  les invariants
$\a_n = \a_n (\cN),
\beta_n= \beta_n(\cN), b_0= b_0(\cN)$ et la fonction $q=q_\cN$ (cf. 2.4). 

  On pose 
$L_2 =
\bigoplus _{k\in
\bZ} R_A(-k)^{l_2(k)}$ et $L_{2, \leq n} = \bigoplus _{k\leq n}
R_A(-k)^{l_2(k)}$ et on appelle $s_n$ la restriction de $s$ \`a $L_{2, \leq n}$ (ou
la matrice correspondante) et $s_{n,t}$ sa valeur au point ferm\'e de $T$. On a alors
les r\'esultats suivants :

\th {Proposition 3.1}.  Soit $n \in \bZ$. \lign
1) On a $\a_n = \rg s_{n,t}$  , de sorte que $\a_n$ est le plus grand entier $\a$ tel
que $s_{n,t}$ ait un $\a$-mineur non nul. \lign
2) L'entier $\beta_n$ est le plus grand entier $\beta$ tel que les $\beta$-mineurs de
$s_{n,t}$ soient sans facteur commun dans $R$. \lign
3) Un entier $n$ est  $\leq b_0$ si et seulement si il v\'erifie:\lign
\indent a) $\a_n= \beta_n$, \lign
\indent b) le sous-$R$-module $F$ de $L_{1,t}$ engendr\'e par les colonnes de la
matrice
$s_{n,t}$ est libre de rang $\a_n$. \lign
En particulier, si $\inf L_2$ est le plus petit entier $n$ tel que $l_2(n)$ soit non nul, on a
$b_0
\geq
\inf L_2$.

\dem Consid\'erons $\s_{N',n} : \cO_\bP(n) \T_k N'_n \fl \cN_t$. Il s'agit de
d\'eterminer le rang de
$\s_{N',n}$ aux points de codimension $\leq 1$ de $\bP^3$. Comme $L_{2,t}\fl N'$ est
surjectif il revient au m\^eme de calculer le  rang de $s'_{n,t}$ o\`u
$s'_{n,t}: (\cL_{2,
\leq n})_t
\fl
\cN_t$ est induite par
$s_{n,t}$. Comme la suite
$(*)$ est localement scind\'ee,  il revient encore au m\^eme de calculer le 
rang de $
\wi{s_{n,t}}$ o\`u $\wi{s_{n,t}}: (\cL_{2, \leq n})_t \fl \cL_{1,t}$  a pour
matrice
$s_{n,t}$ et on en d\'eduit aussit\^ot 1) et 2).

Pour l'assertion 3), supposons d'abord $n \leq b_0$. On a alors $\a_n= \beta_n$ (cf.
1.12). De plus, en vertu  de 1.14,  le
sous-module
$N'_{\leq n}$ de $H^0 \cN_t$  est libre de rang $\a_n$. Mais,
comme $(L_{2, \leq n})_t\fl N'_{\leq n}$  est surjectif et
$\cN_t \fl \cL_{1,t}$ injectif, $N'_{\leq n}$ n'est autre
que $F$.

R\'eciproquement, si on a les conditions a) et b) ci-dessus, on consid\`ere le faisceau $\cL =
\cL_{2,\leq n}$. Comme $F$ est libre de rang $\a_n$, on peut, quitte \`a faire un changement
de base et \`a supprimer  \'eventuellement des colonnes de $s$, supposer que
$\cL$ est dissoci\'e de rang
$\a_n$. La fl\`eche $s_{n,t}$ est alors injective et, comme on a $\a_n= \beta_n$, son
conoyau est sans torsion (cf. [MDP3] II 2.2) et on conclut par 1.14.

\tarte {b) Exemples}

Nous donnons quelques exemples tr\`es simples d'application de la m\'ethode
 ci-dessus. Le lecteur qui souhaiterait avoir plus d'exemples, ou des
exemples plus complexes, ou encore comprendre comment ces exemples sont construits, ou
enfin  quelles sont les  relations de ces exemples avec les sp\'ecialisations dans le
sch\'ema de Hilbert, se reportera  \`a [HMDP3].

\vskip 0.3 cm

Dans tout ce qui suit on suppose que $A$ est une $k$-alg\`ebre qui est un \avd\
d'uniformisante
$a$.

Les exemples ci-dessous sont tous construits de la m\^eme fa\c con : on prend une
suite exacte de $R= k[X,Y,Z,T]$-modules libres gradu\'es : $M_2 \Fl {\s_2} M_1 \Fl
{\s_1} M_0$ telle que le faisceau $\Coker \wi \s_1$ soit localement libre. On obtient
par exemple une telle suite en prenant  trois termes
cons\'ecutifs  dans la r\'esolution libre d'un $R$-module de longueur finie. On
consid\`ere alors $s :
L_2=(M_1
\T_k A)\, \oplus \, (M_2 \T_k A) \fl L_1=(M_0
\T_k A)\, \oplus \, (M_1 \T_k A) $ donn\'e par la matrice
$$s = \pmatrix {\s_1 & 0 \cr a\, I & \s_2 \cr}.$$
On v\'erifie que $\Coker \wi s$ est localement libre, de sorte qu'on peut appliquer 3.1
\`a $\cN= \Im \wi s$ ($L_2 \fl N$ est surjective).

Dans les exemples 3.2 et 3.3 on part de la r\'esolution du $R$-module $k = R/(X,Y,Z,T)$
 donn\'ee par le complexe de Koszul associ\'e \`a la suite r\'eguli\`ere
$U=(X, Y, Z, T)$  : 
$$0 \fl R(-4) \fl R(-3)^4 \Fl {V'} R(-2)^6 \Fl {V} R(-1)^4 \Fl {U} R \fl k \fl 0,$$
$$ V = \pmatrix {Y&Z&T&0&0&0 \cr -X&0&0&Z&T&0 \cr 0&-X&0&-Y&0&T \cr 0&0&-X&0&-Y&-Z
\cr} \quad {\rm et} \quad V' = \pmatrix {0&0&-T&Z\cr
0&T&0&-Y\cr0&-Z&Y&0\cr-T&0&0&X\cr Z&0&-X&0\cr -Y&X&0&0\cr}
$$

\expl {3.2}
On consid\`ere la matrice $s : R_A(-1)^4 \oplus R_A(-2) ^6
\fl R_A \oplus R_A(-1)^4$ avec
$
s = \pmatrix { U & 0\cr
a \,I_4 & V \cr
}$
o\`u $I_4$ est la matrice identit\'e d'ordre $4$.
Au point ferm\'e (i.e., pour $a =0$) il ne reste dans $s_1$ que la matrice ligne $U$ ce
qui montre, en vertu de 3.1,  qu'on a $\a_1= \beta_1 =1$. Comme les colonnes de cette
matrice engendrent l'id\'eal $m= (X, Y, Z, T)$ qui n'est pas libre sur $R$ on a $b_0
=0$. En degr\'e
$2$ on v\'erifie aussit\^ot qu'on a $\a_2= \beta_2 = 4$. On en d\'eduit que la seule valeur
non nulle de la fonction $q$ est $q(2) = 3$ et, en vertu de 2.6, qu'on a une suite
exacte
$$0 \fl \cO_{\bP_A} (-2)^3  \Fl{u} \cN \fl \cJ_{\cC} (2) \fl 0$$
o\`u $\cC$ est une famille de courbes, plate sur $A$ et \`a sp\'ecialit\'e constante. Ces
courbes sont  de degr\'e
$6$ et de genre
$3$. Pour $u$ g\'en\'eral, la courbe g\'en\'erique $C$  de la famille $\cC$ est une
courbe arithm\'etiquement de Cohen-Macaulay (ACM) lisse tandis que la courbe sp\'eciale $C_0$
est une courbe lisse de
bidegr\'e $(2,4)$ trac\'ee sur une quadrique lisse.  Cette courbe est dans la classe de
biliaison de deux droites disjointes et son module de Rao
est concentr\'e en degr\'e $2$.

\expl {3.3}
On consid\`ere  la matrice
$s : R_A (-2)^6 \oplus R_A (-3)^4 \fl R_A (-1)^4 \oplus R_A (-2)^6$ :
$s= \pmatrix {V&0\cr a I_4& V' \cr}$
et le faisceau $\cN$ associ\'e.  On a 
$\a_2= \beta_2= 3$ et $b_0=1$, puis $a_3= \beta_3= 6$ donc la fonction $q$ v\'erifie
$q(2) = 2$ et $q(3) = 3$. La famille minimale de courbes associ\'ee \`a $\cN$ est encore
une famille de courbes de degr\'e
$6$ et genre $3$ comme dans l'exemple 3.2, mais, cette fois-ci, il s'agit d'une famille \`a
postulation constante. Pour un homomorphisme g\'en\'eral, la courbe g\'en\'erique
est une courbe  ACM lisse et  la courbe sp\'eciale est   r\'eunion d'une quartique
plane et de deux droites qui coupent chacune la quartique transversalement en un
point. Cette courbe est dans la classe de biliaison
de deux droites disjointes et son module de Rao a un unique terme non nul en degr\'e
$1$ (et non plus en degr\'e
$2$ comme dans l'exemple 3.2). Pour une
\'etude plus approfondie du sch\'ema de Hilbert
$H_{6,3}$ cf. [AA].

\expl {3.4} Nous donnons ici un exemple avec une partie ``obligatoire'' au
sens de 1.6 c'est-\`a-dire un exemple dans lequel il existe un entier $n\leq b_0$
avec
$q(n)
\neq 0$.

On consid\`ere un module $M$ de longueur finie concentr\'e en degr\'es $0$ et $1$ et
de dimensions respectives $2$ et $7$. On suppose ce module g\'en\'erique et il admet
alors une r\'esolution de la forme suivante  (cf. [MDP2] IV 1.2) 
$$\cdots \fl R(-3)^{34} \Fl {\s_2} R(-1) \oplus R(-2)^{16} \Fl {\s_1} R^2 \fl M \fl 0
\quad {\rm avec} $$
 $$  \s_1\! =\! \pmatrix {\scriptstyle \!X&\scriptstyle \!Y^2&\scriptstyle
Z^2&\scriptstyle \!T^2&\scriptstyle \! YZ&\scriptstyle \!\scriptstyle \!YT&\scriptstyle
ZT&\scriptstyle \! 0&\scriptstyle \!0&\scriptstyle \!0&\scriptstyle \!0&\scriptstyle
0&\scriptstyle \!0&\scriptstyle \! 0&\scriptstyle \!0&\scriptstyle \!0&\scriptstyle \!0\cr
\scriptstyle \!-Y&\scriptstyle \!0&\scriptstyle \! 0&\scriptstyle \!0&\scriptstyle \!0&\scriptstyle
0&\scriptstyle \!0&\scriptstyle \! X^2&\scriptstyle \!Y^2&\scriptstyle \!Z^2&\scriptstyle
T^2&\scriptstyle \! XY&\scriptstyle \!XZ&\scriptstyle \!XT&\scriptstyle \!YZ&\scriptstyle
YT&\scriptstyle \!ZT \cr}$$ On consid\`ere alors 
$ s:  R_A(-1) \oplus R_A(-2)^{16} \oplus  R_A(-3)^{34}  \Fl {s} R_A^2 \oplus
R_A(-1)
\oplus R_A(-2)^{16}  $
$$ \hbox {donn\'e par} \qquad s= \pmatrix {\s_1 & 0 \cr a\, I_{17} & \s_2 \cr}.$$
La matrice $s$ admet une seule colonne en degr\'e $\leq 1$, dont la
transpos\'ee est la ligne $(X, \,-Y,\,  -a, \,  0, \cdots,\, 0)$. On en d\'eduit
aussit\^ot en utilisant 3.1 qu'on a $\a_1= \beta_1 = 1$ et $b_0 \geq 1$  et donc, par
d\'efinition de la fonction $q$, $q(1) =1$. On a aussi $q(3) = 15$ et la famille minimale est
en degr\'e $120$ et genre $1001$ avec une courbe g\'en\'erique ACM.

\titre {R\'ef\'erences bibliographiques}

[AA] A\"\i t-Amrane S., Sur le sch\'ema de Hilbert $H_{d, (d-3)(d-4)/2}$, en
pr\'eparation.

[AG] Hartshorne R., Algebraic geometry, Graduate texts in Mathematics 52, Springer
Verlag, 1977.

[EGA] Grothendieck A. et Dieudonn\'e J., \'El\'ements de g\'eom\'etrie alg\'ebrique III,
Publ. Math. IHES 11, 1961 et IV,  Publ. Math. IHES 28, 1966.

[H] Hartshorne R., Coherent functors, \`a para\^\i tre, Advances in Math.

[HMDP1] Hartshorne R.,  Martin-Deschamps M. et  Perrin D., Un th\'eor\`eme de Rao pour
les familles de courbes gauches, en pr\'eparation.

[HMDP3] Hartshorne R.,  Martin-Deschamps M. et  Perrin D., Triades et familles de courbes
gauches, en pr\'eparation.

[K] Kleiman S.L., Geometry on Grassmannians and applications to splitting bundles and
smoothing cycles, Publ. Math. IHES 36, 1969, 281-297.

[LR] Lazarsfeld R. et Rao A. P., Linkage of general curves of large degree, Lecture
notes 997, Springer Verlag, 1983, 267-289.

[MDP 1]  Martin-Deschamps M. et  Perrin D., Sur la classification des
courbes gauches, Ast\'erisque, Vol. 184-185, 1990.

[MDP 2] Martin-Deschamps M. et  Perrin D., Courbes gauches et Modules de Rao,
J. reine angew. Math. 439, 1993, 103-145.

[MDP3]  Martin-Deschamps M. et  Perrin D., Construction de courbes lisses :
un th\'eor\`eme \`a la Bertini, rapport de recherche du LMENS 92-22, 1992.

[MDP4]  Martin-Deschamps M. et  Perrin D., Quand un morphisme de fibr\'es
d\'eg\'en\`ere-t-il le long d'une courbe lisse, \`a para\^\i tre dans les actes des
conf\'erences Europroj de Catane et Barcelone, Marcel Dekker, Inc.

 [R] Rao A.P., Liaison equivalence classes, Math. Ann. 258, 1981, 169-173.

[S1] Schlesinger E., Vms fundamentals, unpublished manuscript, Berkeley, 1992.

 \bye